\documentclass[a4paper,12pt,reqno]{amsart}
\usepackage{pb-diagram} 
\usepackage{amsmath,amssymb}

\usepackage{amsmath}
\usepackage{enumerate}
\usepackage{amsthm}
\usepackage{hyperref}
\usepackage[active]{srcltx}
\usepackage{pb-diagram}

\newtheorem{theorem}{Theorem}
\newtheorem{proposition}[theorem]{Proposition}
\newtheorem{corollary}[theorem]{Corollary}
\newtheorem{remark}{Remark}
\newtheorem{lemma}[theorem]{Lemma}
\newfont{\bb}{msbm10 at 12pt}

\DeclareMathAlphabet{\doba}{U}{msb}{m}{n}   
 
\def\what{\widehat}
\def\E{\mathcal{E}}
     
\def\SOm{\mathbb{S}(\Omega)}

\def\SSi{\mathbf{S}}     
\def\Id{{\mathop{\rm Id}}} 
     
\def\SSi{\mathbf{S}}     
\def\Id{{\mathop{\rm Id}}}     
\def\<{\langle}     
\def\>{\rangle}     
\def\div{{\rm div}}
   
\newcommand{\RM}{\mathbb{R}^{3,1}}
\newcommand{\SM}{\mathbb{S}\mathbb{R}^{3,1}}

\long\def\komment#1{} 

\begin{document}


\title[]{The Dirac operator on untrapped surfaces}     
\author{Simon Raulot}     
\address
{Laboratoire de Math\'ematiques R. Salem\\
UMR $6085$ CNRS-Universit\'e de Rouen\\
Avenue de l'Universit\'e, BP.$12$\\
Technop\^ole du Madrillet\\
$76801$ Saint-\'Etienne-du-Rouvray, France}
\email{simon.raulot@univ-rouen.fr}
\date{\today}

\maketitle

\begin{abstract}
We establish a sharp extrinsic lower bound for the first eigenvalue of the Dirac operator of an untrapped surface in initial data sets without apparent 
horizon in terms of the norm of its mean curvature vector. The equality case leads to rigidity results for the constraint equations with spherical boundary
as well as uniqueness results for constant mean curvature surfaces in Minkowski space.
\end{abstract}


\section{Introduction}


In \cite{friedrich}, T. Friedrich proves that on a closed connected $n$-dimensional Riemannian spin manifold $(M^n,g)$, the first eigenvalue of the Dirac operator 
satisfies 
\begin{eqnarray*}
\lambda_1^2\geq\frac{n}{4(n-1)}\inf_{M} R,
\end{eqnarray*}

where $R$ is the scalar curvature of $(M^n,g)$. Then, for $n\geq 3$, using the conformal behavior of the Dirac operator, O. Hijazi proves in \cite{hijazi:86} that, up to a 
dimensional constant, $\lambda^2_1$ is at least equal to the first eigenvalue of the Yamabe operator. For the case of $2$-spheres, C. B\"ar
\cite{baer:92} (see also \cite{hijazi:91}) shows that 
\begin{eqnarray*}
 \lambda_1^2\geq\frac{4\pi}{{\rm Area}(M)}
\end{eqnarray*}

and the equality holds only for the round metric. More recently, the spectrum of the Dirac operator has been studied from an extrinsic point of view. Indeed, 
in \cite{hijazi.montiel.zhang:1} and \cite{hijazi.montiel.roldan:02}, the authors prove that if $\Sigma$ is a $n$-dimensional compact hypersurface which bounds
a compact domain $\Omega$ in an $(n+1)$-dimensional Riemannian spin manifold and if the scalar curvature of $\Omega$ is such that $R\geq -n(n+1)\alpha^2$, 
with $\alpha\in\mathbb{R}$, then the first positive eigenvalue $\lambda_1(D^\Sigma)$ of the Dirac operator on $\Sigma$, satisfies
\begin{eqnarray}\label{extri}
 \lambda_1(D^\Sigma)\geq\frac{1}{2}\inf_{\Sigma}\Big(\sqrt{H^2_\Sigma-n^2\alpha^2}\Big)
\end{eqnarray}
where $H_\Sigma$ is the mean curvature of $\Sigma$ in $\Omega$ satisfying $H_\Sigma\geq n|\alpha|$. Note that here the mean curvature is the trace of the second fundamental 
form and not the normalized trace as in \cite{hijazi.montiel.zhang:1} and \cite{hijazi.montiel.roldan:02}. Several geometric applications, like Alexandrov type theorems
or rigidity results for compact manifolds with boundary, are derived from these estimates (see also \cite{rollmops}).\\

In this paper, we generalize the inequality (\ref{extri}) for untrapped surfaces in initial data sets of the Einstein equation. From this estimate, we
deduce several rigidity results for the constraint equation with spherical boundary as well as uniqueness results for constant mean curvature surfaces in 
the Minkowski spacetime.


\subsection{The main estimate}


Here we introduce the setting and give the precise statement of our main result.

\medskip

Let $(\E^4,g_\E)$ be a spacetime satisfying the Einstein field equations, that is $(\E^4,g_\E)$ is a $4$-dimensional time oriented Lorentzian manifold such that
\begin{eqnarray*}
Ric_\E-\frac{1}{2} R_\E g_\E=\mathcal{T},
\end{eqnarray*}

where $R_\E$ (resp. $Ric_\E$) denotes the scalar curvature (resp. the Ricci curvature) of ($\E,g_\E)$ and $\mathcal{T}$ is the 
energy-momentum tensor of the spacetime. We also assume that $(\E^4,g_\E)$ satisfies the dominant energy condition that is its energy-momentum tensor
$\mathcal{T}$ has the property that the vector field dual to the one form $-\mathcal{T}(\eta,.)$ is a future direct causal vector of $T\E$, for every future direct causal
vector $\eta\in\Gamma(T\E)$. 

\medskip

Let $M^3$ be an immersed spacelike hypersurface of $(\E^4,g_\E)$ with induced Riemannian metric $g$. Assume that $T$ is the future directed timelike normal vector 
to $M$ and denote by $K$ the associated second fundamental form defined by $K(X,Y)=g_\E(\nabla^\E_X T,Y)$, for all $X,Y\in\Gamma(TM)$. Here $\nabla^\E$ denotes the Levi-Civita connection of the spacetime. Then the Gauss, Codazzi and Einstein equations provide {\it constraint equations} 
on $M$ given by  
$$
\begin{array}{rll}
\frac{1}{2}\big(R-|K|^2_M+({\rm Tr}_M(K))^2\big) & = & \mu\\
-\div_M\big(K-{\rm Tr}_M(K)g\big) & = & J
\end{array}
$$
where $R$ is the scalar curvature of $(M^3,g)$, $|K|^2$ and ${\rm Tr}(K)$ denote the squared norm and the trace of $K$ with respect to $g$ and $\mu$ and $J$ 
are respectively the energy and the momentum densities of the matter fields given by some components of the energy-momentum tensor $\mathcal{T}$. 

A consequence of the dominant energy condition for the spacetime is that the inequality
\begin{eqnarray*}
 \mu\geq |J|
\end{eqnarray*}
holds on $M$. A triplet $(M^3,g,K)$ which satisfies the dominant energy condition is called an {\it initial data set}.

\medskip

Now we consider a spacelike $2$-surface $\Sigma$ in $(\E^4,g_\E)$ which bounds a compact domain $\Omega$ in the initial data set $(M^3,g,K)$. Denote by
$N$ the inward unit vector field normal to $\Sigma$ in $\Omega$ so that $T$ and $N$ span the normal bundle $\mathcal{N}\Sigma$ of $\Sigma$ in $\E$. 
The second fundamental form of $\Sigma$ in $\E$ is defined by
\begin{eqnarray*}
 II(X,Y)=-A(X,Y) N-K(X,Y) T,
\end{eqnarray*}
for all $X,Y\in\Gamma(T\Sigma)$ and where $A(X,Y):=-g(\nabla_XN,Y)$ denotes the second fundamental form of $\Sigma$ in $\Omega$. Here $\nabla$ denotes the 
Levi-Civita connection of the Riemannian metric $g$ on $M$. The trace of the second fundamental form $II$ along the surface $\Sigma$ defines the 
mean curvature vector field $\mathcal{H}$ of $\Sigma$ in $\E$ and it is given by:
\begin{eqnarray*}
 \mathcal{H}=-H_\Sigma N-{\rm Tr}_{\Sigma}(K)T.
\end{eqnarray*}
Its norm can be expressed as: 
$$
|\mathcal{H}|^2= H_\Sigma^2-\mathrm{Tr}_{\Sigma}(K)^2=\theta_+\theta_-
$$
where $H_\Sigma:={\rm Tr}_\Sigma(A)$ is the mean curvature of $\Sigma$ in $\Omega$ and $\theta_\pm=H_\Sigma\pm\mathrm{Tr}_{\Sigma}(K)$ are its null expansions. 
The spacelike surfaces with $\theta_+<0$ (or $\theta_-<0$) are referred to as future (or past) trapped surfaces. A surface with $\theta_+=0$ or 
$\theta_-=0$ is called an {\it apparent horizon}. Finally, if $\Sigma$ is such that $\theta_+>0$ and $\theta_->0$ (that is $H_\Sigma>|{\rm Tr}_\Sigma(K)|$), 
it will be referred to as an untrapped surface. In particular, an untrapped surface has a spacelike mean curvature vector. 

\medskip

We now give the precise statement of our main result:
\begin{theorem}\label{spacetime}
Let $\Omega$ be a compact domain with smooth boundary $\Sigma$ in an initial data set $(M^3,g,K)$ and assume that $\Omega$ has no apparent horizon in its interior. 
If $\Sigma:=\partial\Omega$ is an untrapped surface in $(\E^4,g_\E)$, then the first positive eigenvalue $\lambda_1(D^\Sigma)$ of the Dirac 
operator $D^\Sigma$ of $\Sigma$ satisfies
\begin{eqnarray}\label{eigenspacetime}
\lambda_1(D^\Sigma)\geq\frac{1}{2}\inf_{\Sigma}\Big(\sqrt{H_\Sigma^2-\mathrm{Tr}_{\Sigma}(K)^2}\Big).
\end{eqnarray}

Moreover, if equality occurs then $(\Omega^3,g)$ can be locally embedded into the Minkowski flat spacetime as a spacelike hypersurface with $K_{ij}$ as 
second fundamental form.
\end{theorem}

Using an approach proposed in \cite{miaoshitam}, we prove (see Theorem \ref{Mink}) that this estimate holds in the Minkowski spacetime without any assumption 
on the existence of apparent horizons. Moreover, in this context, we show that equality holds in (\ref{eigenspacetime}) only for round spheres. This rigidity
result combined with a new upper bound for $\lambda_1(D^\Sigma)$ (see Proposition \ref{upper}) applies to obtain an Alexandrov type theorem (see Corollary 
\ref{genalex}) as well as a unified proof of the Euclidean and Hyperbolic versions of the Alexandrov theorem for surfaces (see Corollary \ref{alex}).  

\medskip

Let us briefly discuss the proof of Theorem \ref{spacetime}. Here we use arguments developed by Schoen-Yau \cite{sy1} in their proof of the positive mass theorem 
for general initial data sets and by Liu-Yau \cite{liuyau} in their proof of the positivity of a quasi-local mass.

\medskip

Indeed in a first part, we prove a generalization of the inequality (\ref{extri}) with $\alpha=0$ for Riemannian manifolds on which a non negativity condition 
holds on the scalar curvature. More precisely, we have
\begin{theorem}\label{newextrin}
Let $\Sigma$ be a compact hypersurface in a Riemannian spin manifold $M$ bounding a 
compact domain $\Omega$. Assume that the scalar curvature $R$ of the domain $\Omega$ and the mean curvature $H_\Sigma$ of $\Sigma$ satisfy
\begin{eqnarray*}
R\geq 2|X|^2+2\div(X)\quad and \quad H_\Sigma\geq -g(X,N).
\end{eqnarray*}  

The lowest nonnegative eigenvalue $\lambda_1(D^\Sigma)$ of the Dirac operator $D^\Sigma$ satisfies
\begin{eqnarray}\label{ne1}
\lambda_1(D^\Sigma)\geq\frac{1}{2}\underset{\Sigma}{\inf}\big(H_\Sigma+g(X,N)\big).
\end{eqnarray}

Moreover equality occurs if and only if $\Omega$ carries a parallel spinor, $\Sigma$ has constant mean curvature, $X\equiv 0$ on $\Omega$ and the eigenspace 
corresponding to $\lambda_1(D^\Sigma)$ consists of the restrictions of parallel spinors on the domain $\Omega$.
\end{theorem}

\medskip Then, starting from an untrapped surface $\Sigma$ which bounds a compact domain $\Omega$ in an initial data set $(M^3,g,K)$, we use a solution of the Jang 
equation with Dirichlet boundary condition, to get a new metric $\what{g}$ on $\Omega$ for which  Theorem \ref{newextrin} applies. The conclusion follows from
boundary calculations of the mean curvature of $\Sigma$ in $(\Omega,\what{g})$ due to Yau \cite{yau1}.

\begin{remark}
If $\Sigma$ is the boundary of a domain $\Omega$ in an initial data set $(M^3,g,\alpha g)$ with $\alpha\in\{0,1\}$, the dominant 
energy condition reads as $R_M\geq -6\alpha^2$ and the untrapped condition gives $H_\Sigma>2\alpha$. We can then apply  Theorem \ref{spacetime} 
to get
\begin{eqnarray*}
\lambda_1(D^\Sigma)\geq \frac{1}{2}\inf_{\Sigma}\Big(\sqrt{H_\Sigma^2-4\alpha^2}\Big)
\end{eqnarray*}

which is precisely the eigenvalue estimates proved in \cite{hijazi.montiel.zhang:1} for $\alpha=0$ and in \cite{hijazi.montiel.roldan:02} for $\alpha=1$.
\end{remark}

{\bf Acknowledgements:} I would like to thank Oussama Hijazi and Emmanuel Humbert for their remarks and suggestions, as well as their support.


\section{A new extrinsic lower bound in the Riemannian setting}



\subsection{Preliminaries}


In this section, we consider a compact Riemannian spin manifold $(\Omega^n,g)$ with smooth boundary $\Sigma$ (possibly with several connected components) and 
denote by $\nabla$ the Levi-Civita connection on the tangent bundle $T\Omega$. We fix a spin structure on $\Omega$ and denote by $\SOm$ the corresponding 
spinor bundle. This complex vector bundle is naturally endowed with a left Clifford multiplication
\begin{eqnarray*}
\gamma:\mathbb{C}\ell(\Omega)\longrightarrow \hbox{End}(\SOm)
\end{eqnarray*}
which is a fiber preserving algebra morphism. Then $\SOm$ becomes a bundle of complex left modules over the Clifford bundle $\mathbb{C}\ell(\Omega)$ of 
$\Omega$. The spinor bundle $\SOm$ is equipped with a natural Hermitian metric, denoted by $\langle\; ,\;\rangle$, and with the spinorial Levi-Civita connection
$\nabla$ acting on spinor fields which are compatible with the Clifford multiplication $\gamma$. The Dirac operator
$D$ on $\SOm$ is  locally given by: 
$$
D=\sum_{i=1}^{n}\gamma(e_i)\nabla_{e_i},
$$
where $\{e_1,\dots,e_{n}\}$ is a local orthonormal frame of $T\Omega$ and it defines a first order elliptic differential operator. On the boundary $\Sigma$, we denote also by $g$ the induced Riemannian metric 
and by $\nabla^\Sigma$ its Levi-Civita connection. Recall that since $\Sigma$ is oriented, it admits a spin structure induced from the one on $\Omega$. Hence we 
have that the restriction
\begin{eqnarray*}
\SSi:=\SOm_{|\Sigma}
\end{eqnarray*}
is a left module over $\mathbb{C}\ell(\Sigma)$ with Clifford multiplication given by
\begin{eqnarray*}
\gamma^\SSi(X)\psi=\gamma(X)\gamma(N)\psi
\end{eqnarray*}
for every $\psi\in\Gamma(\SSi)$ and $X\in\Gamma(T\Sigma)$. Here $N$ denotes the unit inward vector field normal to $\Sigma$. Now the Riemannian Gauss formula 
implies that the spin connection $\nabla^\SSi$ defined on $\SSi$ by 
\begin{equation}\label{spin-gaus}
\nabla^\SSi_X\psi=\nabla_X\psi-\frac{1}{2}\gamma^\SSi(A X)\psi
\end{equation}
for $X\in\Gamma(T\Sigma)$ and the Hermitian metric $\langle\; ,\;\rangle$ induced from that of $\SOm$ are compatible with the Clifford multiplication 
$\gamma^\SSi$. Here $A(X):=-\nabla_XN$ denotes the Weingarten map of $(\Sigma,g)$ in $(\Omega,g)$. The extrinsic Dirac 
operator of $\Sigma$ acting on $\SSi$ is then defined by $D^\SSi=\gamma^\SSi\circ\nabla^\SSi$. It is a well known fact that $D^\SSi$ is a first order 
elliptic differential operator which is formally $L^2$-selfadjoint. By (\ref{spin-gaus}), it is straightforward to compute that for all 
$\psi\in \Gamma(\SSi)$:
\begin{equation}\label{dira-extr2}
D^\SSi\psi=\frac{1}{2}H_\Sigma\psi-\gamma(N)D\psi-\nabla_N\psi.
\end{equation}
Here $H_\Sigma:={\rm Tr}(A)$ denotes the mean curvature of $\Sigma$ in $\Omega$ for the metric $g$. In our conventions, the unit $2$-sphere in the $3$-dimensional 
Euclidean space has constant mean curvature equal to $2$. It is important to note that since $\Sigma$ is spin, we can also define an intrinsic spinor bundle
$\mathbb{S}(\Sigma)$ over $\Sigma$ and an associated Dirac operator $D^\Sigma$ related with the extrinsic spinor bundle by the following identification (see 
\cite{hijazi.montiel.zhang:1} for more details):
$$(\SSi,D^\SSi)\equiv
\left\lbrace
\begin{array}{ll}
(\mathbb{S}(\Sigma),D^\Sigma) & \quad \text{if n is even}\\
(\mathbb{S}(\Sigma)\oplus\mathbb{S}(\Sigma),D^\Sigma\oplus -D^\Sigma) &  \quad \text{if n is odd}.
\end{array}\right.$$

Recall that in the context of Riemannian spin manifold with boundary, the classical Schr\"odinger-Lichnerowicz formula asserts that:
\begin{eqnarray*}
D^2=\nabla^*\nabla+\frac{R}{4}
\end{eqnarray*}

and it gives the integral formula:
\begin{eqnarray}\label{lichneint}
\int_{\Omega}\big(|\nabla\varphi|^2+\frac{R}{4}|\varphi|^2-|D\varphi|^2\big)dv=\int_{\Sigma}\<D^\SSi\varphi-\frac{1}{2}H_\Sigma\varphi,\varphi\>ds
\end{eqnarray}

for all $\varphi\in\Gamma\big(\mathbb{S}(\Omega)\big)$. Here $\nabla^*$ denotes the $L^2$-adjoint of $\nabla$ and $R$ is the scalar curvature of $(\Omega,g)$.

\vspace{0.2cm}

{}From now, we assume that there exits a smooth vector field $X\in\Gamma(T\Omega)$ such that:
\begin{eqnarray}\label{scalarcond}
R\geq 2|X|^2+2\div(X)
\end{eqnarray}

where $|X|^2=g(X,X)$ and the divergence of a vector field $X=\sum_{j=1}^nX^je_j\in\Gamma(T\Omega)$ is locally defined by
\begin{eqnarray*}
\div(X)=-\sum_{j=1}^n\nabla_{e_i}X_i.
\end{eqnarray*}

We first generalize an argument of Liu-Yau \cite{liuyau} and prove a modified version of Formula (\ref{lichneint}) appropriate to the assumption 
(\ref{scalarcond}). More precisely, we have:
\begin{proposition}\label{modlich}
Let $(\Omega^n,g)$ a compact Riemannian spin manifold with boundary such that there exists a smooth vector field $X\in\Gamma(T\Omega)$ satisfying 
(\ref{scalarcond}), then 
\begin{eqnarray}\label{modlich1}
\int_{\Omega}\big(\frac{1}{2}|\nabla\varphi|^2-|D\varphi|^2\big)dv\leq\int_{\Sigma}\<D^\SSi\varphi-\frac{1}{2}\big(H_\Sigma+g(X,N)\big)\varphi,\varphi\>ds
\end{eqnarray}
for all $\varphi\in\Gamma(\SOm)$. Moreover, equality occurs if and only if the spinor field $\varphi$ satisfies 
\begin{eqnarray*}
\nabla_Y\varphi=-g(X,Y)\varphi
\end{eqnarray*}
for all $Y\in\Gamma(T\Omega)$.
\end{proposition}

{\it Proof:}
First note that since 
\begin{eqnarray*}
 \div(|\varphi|^2 X)=-X(|\varphi|^2)+|\varphi|^2\div(X),
\end{eqnarray*}

the Stokes formula gives
\begin{eqnarray*}
\int_\Omega\frac{R}{4}|\varphi|^2dv & = & \int_\Omega\big(\frac{R}{4}-\frac{1}{2}\div(X)\big)|\varphi|^2dv+\frac{1}{2}\int_\Omega\div(X)|\varphi|^2dv\\
& = & \frac{1}{4}\int_\Omega\big(R-2\div(X)\big)|\varphi|^2dv+\frac{1}{2}\int_\Omega X(|\varphi|^2)dv+\frac{1}{2}\int_\Sigma g(X,N)|\varphi|^2ds.
\end{eqnarray*}

Inserting this identity in the integral form (\ref{lichneint}) of the Schr\"odinger-Lichnerowicz formula leads to:
\begin{eqnarray*}
\int_{\Omega}\big(|\nabla\varphi|^2+\frac{1}{4}\big(R-2\div(X)\big)|\varphi|^2+\frac{1}{2}X(|\varphi|^2)-|D\varphi|^2\big)dv=
\int_{\Sigma}\<D^\SSi\varphi-\frac{1}{2}\big(H_\Sigma+g(X,N)\big)\varphi\>ds
\end{eqnarray*}

and using (\ref{scalarcond}), we conclude that:
\begin{eqnarray}\label{schrint}
\int_{\Omega}\big(|\nabla\varphi|^2+\frac{1}{2}|X|^2|\varphi|^2+\frac{1}{2}X(|\varphi|^2)-|D\varphi|^2\big)dv\leq
\int_{\Sigma}\<D^\SSi\varphi-\frac{1}{2}\big(H_\Sigma+g(X,N)\big)\varphi\>ds.
\end{eqnarray}

If we let $\widetilde{\nabla}_Y\varphi:=\nabla_Y\varphi+g(X,Y)\varphi$, it is straightforward to compute
\begin{eqnarray*}
|\widetilde{\nabla}\varphi|^2=|\nabla\varphi|^2+|X|^2|\varphi|^2+2{\rm Re}\<\nabla_X\varphi,\varphi\> 
\end{eqnarray*}

and since $2{\rm Re}\<\nabla_X\varphi,\varphi\>=X(|\varphi|^2)$ we get:
\begin{eqnarray*}
\frac{1}{2}X(|\varphi|^2)\geq -\frac{1}{2}|\nabla\varphi|^2-\frac{1}{2}|X|^2|\varphi|^2
\end{eqnarray*}

with equality if and only if $\widetilde{\nabla}_Y\varphi=0$ for all $Y\in\Gamma(T\Omega)$. Combining this last inequality with (\ref{schrint}) leads to the 
result.
\hfill$\square$\\


\subsection{The Riemannian estimate}


Now following \cite{hijazi.montiel.zhang:2}, we use the fact that the Dirac operator of the domain $\Omega$ defines an isomorphism between suitable spaces to
extend a spinor field on $\Sigma$ harmonically on $\Omega$. More precisely, we consider the $(MIT)$ bag boundary condition which defines an elliptic condition
for the Dirac operator of $\Omega$. It is given by the pointwise orthogonal projection
$$\begin{array}{lccl}
B^{\pm}: & {\rm L}^2({\bf S}) & \longrightarrow & {\rm L}^2(V^{\pm})\\ 
 & \varphi & \longmapsto & \frac{1}{2}(\Id\pm i\gamma(N)\varphi)
\end{array}$$ 

where $V^{\pm}$ is the eigensubbundles of $\SSi$ associated with the eigenvalues $\pm 1$ of the involution 
$$\begin{array}{lccl}
\mathcal{I}: & {\bf S} & \longrightarrow & {\bf S}\\ 
 & \varphi & \longmapsto & i\gamma(N)\varphi.
\end{array}$$ 

We then recall:
\begin{lemma}\label{existenceMIT}(\cite{hijazi.montiel.zhang:2})
Let $(\Omega^n,g)$ be a compact Riemannian spin manifold with smooth boundary $\Sigma$. Then for all $\Psi\in\Gamma({\bf S})$, the boundary problem 
$$\left\lbrace
\begin{array}{ll}
D\phi = 0 & \qquad\text{on}\,\,\Omega\\
B^{\pm}\phi_{|\Sigma}= B^{\pm}\Psi& \qquad\text{along}\,\,\Sigma
\end{array}
\right.$$

has a unique smooth solution $\phi\in\Gamma\big(\mathbb{S}(\Omega)\big)$.
\end{lemma}

\vspace{0.2cm}

{\it Proof of Theorem \ref{newextrin}:}
Let $\Phi_1\in\Gamma({\bf S})$ be a smooth eigenspinor for the Dirac operator $D^\SSi$ associated with the eigenvalue $\lambda_1(D^\SSi)$. Lemma 
\ref{existenceMIT} ensures the existence of a smooth spinor field $\Psi\in\Gamma\big(\mathbb{S}(\Omega)\big)$ satisfying the boundary value problem:
$$\left\lbrace
\begin{array}{ll}
D\Psi = 0 & \qquad\text{on}\,\,\Omega\\
B^{\pm}\Psi_{|\Sigma}= B^{\pm}\Phi_1& \qquad\text{along}\,\,\Sigma.
\end{array}
\right.$$

Using Formula (\ref{modlich1}) leads to:
\begin{eqnarray}\label{truc}
0\leq\frac{1}{2}\int_{\Omega}|\nabla\varphi|^2dv & \leq & \int_{\Sigma}\<D^\SSi\varphi-\frac{1}{2}\big(H_\Sigma+g(X,N)\big)\varphi,\varphi\>ds\\
& \leq & \Big(\lambda_1(D^\SSi)-\frac{1}{2}\underset{\Sigma}{\inf}\big(H_\Sigma+g(X,N)\big)\Big)\int_{\Sigma}|\varphi|^2ds\nonumber
\end{eqnarray}

which gives the first part of the result. The last inequality is easily derived from the boundary condition (see \cite{hijazi.montiel.zhang:2} or 
\cite{rollmops} for more details). Moreover equality occurs in this inequality if and only if $\Psi_{|\Sigma}=\Phi_1$. Assume now that equality is achieved, 
then because of (\ref{truc}) the spinor $\Psi\in\Gamma\big(\mathbb{S}(\Omega)\big)$ is a parallel spinor such that:
$$\left\lbrace
\begin{array}{ll}
\nabla\Psi = 0 & \qquad\text{on}\,\,\Omega\\
\Psi_{|\Sigma}= \Phi_1& \qquad\text{along}\,\,\Sigma.
\end{array}
\right.$$

On the other hand, equality is also achieved in (\ref{modlich1}) then:
\begin{eqnarray*}
\nabla_Y\Psi=-g(X,Y)\Psi
\end{eqnarray*}

for all $Y\in\Gamma(T\Omega)$. Since $\Psi$ has no zeros, we finally have $X\equiv 0$ on $\Omega$. The end of the proof is then similar to Theorem $5$ in 
\cite{hijazi.montiel.zhang:1}.
\hfill$\square$

\begin{remark}
Theorem~\ref{newextrin} is a generalization of Theorem $5$ in \cite{hijazi.montiel.zhang:1}. Indeed, for $X=0$ we get precisely their result.
\end{remark}

To conclude this section, it should be pointed out that the Theorem \ref{newextrin} can be viewed as a corollary of a more general result:
\begin{theorem}\label{FunctEst}
Assume that the assumptions of Theorem \ref{newextrin} are fulfilled. Let $H_0\in C^\infty(\Sigma)$ be a smooth nonnegative function on $\Omega$ and 
$\Phi\in\Gamma({\bf S})$ a smooth spinor field such that
\begin{eqnarray}\label{diracequation}
D^\SSi\Phi=H_0\Phi
\end{eqnarray}

then
\begin{eqnarray*}
\frac{1}{2}\int_{\Sigma}\big(H_\Sigma+g(X,N)\big)|\Psi|^2ds\leq\int_{\Sigma}H_0|\Psi|^2ds
\end{eqnarray*}

where $\Psi\in\Gamma\big({\mathbb S}(\Omega)\big)$ is the unique harmonic extension of $\Phi$ for the (MIT) condition. Moreover, equality occurs if and only 
if $\Omega$ carries a parallel spinor, $\Sigma$ has mean curvature $H=2H_0$, $X\equiv 0$ on $\Omega$ and the solutions of (\ref{diracequation}) consist of 
restrictions of parallel spinors on the domain $\Omega$.
\end{theorem}

Theorem \ref{newextrin} follows for $\Phi=\Phi_1$, an eigenspinor of the Dirac operator $D^\SSi$ associated with the first nonnegative eigenvalue 
$\lambda_1(D^\SSi)$ (that is $H_0=\lambda_1$).


\section{Untrapped surfaces in initial data sets}\label{casgeneral}



\subsection{The Jang Equation}


In this section, we recall some well-known facts on the Jang equation. For more details on this subject, we refer to \cite{sy1}, \cite{yau1} or \cite{AEM}.
This equation first appears in \cite{jang} aiming to prove the positive mass theorem using the inverse mean curvature flow but without success. However,
as shown by Schoen and Yau \cite{sy1}, this equation can be used to reduce the proof of the general positive mass theorem to the case of time-symmetric
initial data sets (that is $K_{ij}=0)$ previously obtained by the same authors in \cite{symac}. More recently, Yau and Liu \cite{liuyau} defines a quasi-local mass, generalizing the Brown-York 
quasi-local mass, and prove its positivity using the Jang equation.

\medskip

The problem can be stated as follow: let $(M^3,g,K)$ be an initial data set for the Einstein equation and consider the four dimensional manifold 
$M\times\mathbb{R}$ equipped with the {\it Riemannian} metric $\<\,,\,\>:=g\oplus dt^2$. The problem is to find a smooth function $u:M\rightarrow\mathbb{R}$   
such that the hypersurface $\widehat{M}$ of $M\times\mathbb{R}$ obtained by taking the graph of $u$ over $M$ satisfies the equation
\begin{eqnarray*}
H_{\what{M}}={\rm Tr}_{\what{M}}(K)
\end{eqnarray*}
where $H_{\what{M}}$ denotes the mean curvature of $\what{M}$ in $(M\times\mathbb{R},\<\,,\,\>)$ and ${\rm Tr}_{\what{M}}(\,.\,)$ is the trace on $\what{M}$ 
with respect to the induced metric. This geometric problem is equivalent to solve the non-linear second order elliptic equation
\begin{eqnarray}\label{jang}
\sum_{i,j=1}^3\Big(g^{ij}-\frac{u^iu^j}{1+|\nabla u|^2}\Big)\Big(\frac{(\nabla^2u)_{ij}}{\sqrt{1+|\nabla f|^2}}-K_{ij}\Big)=0 
\end{eqnarray}
where $\nabla$ (resp. $\nabla^2$) denotes the Levi-Civita connection (resp. the Hessian) of the metric $g$, $u^i=g^{ij}u_j$ and $u_j=e_j(u)$. Note that the 
metric induced by $\<\,,\,\>$ on $\what{M}$ is
\begin{eqnarray*}
 \what{g}_{ij}=g_{ij}+u_iu_j
\end{eqnarray*}
and can be viewed as a deformation of the metric $g$ on $M$. In the following, we adopt the convention that $M$ and $\what{M}$ denote respectively the Riemannian 
manifolds $(M,g)$ and $(M,\what{g})$. Analogously, if $\nabla$ denotes the Levi-Civita connection for $M$, then $\what{\nabla}$ denotes the Levi-Civita 
connection on $\what{M}$ and so on. Assuming now that the initial data set $(M^3,g,K)$ comes from a spacetime satisfying the dominant energy condition,
Schoen and Yau proved that the following relation holds on $\what{M}$:
\begin{eqnarray}\label{schoenyau}
 0\leq 2(\mu-|J|)\leq \what{R}-2|X|_{\what{g}}^2-2\what{\div}(X)
\end{eqnarray}
where 
\begin{eqnarray}\label{defX}
X=\omega-\what{\nabla}\log(f),
\end{eqnarray}
$\omega$ is the tangent part of the vector field dual to $-K(\,.\,,\what{\nu})$, $f=-\<\partial_t,\what{\nu}\>$ and $\what{\nu}$ denotes the unit normal 
vector field to $\what{M}$ in $M\times\mathbb{R}$. All the quantities $K_{ij}$, $\mu$ and $J$ are defined on $M\times\mathbb{R}$ by parallel transport along 
the $\mathbb{R}$-factor. Moreover, equality occurs in (\ref{schoenyau}) if and only if $\mu=|J|$ and the second fundamental form of $\what{M}$ in 
$M\times\mathbb{R}$ is $K$. 

\vspace{0.2cm}

It is important to note here that in  Theorem \ref{spacetime}, we assume that there is no apparent horizon in the interior of $\Omega$ to ensure the existence of a 
global solution to the Jang equation. Indeed, if $\widetilde{\Sigma}$ is an apparent horizon then a solution of (\ref{jang}) blows up at $\widetilde{\Sigma}$ and the resulting graph 
$\what{M}$ is asymptotic to a cylinder $\widetilde{\Sigma}\times \mathbb{R}$. To overcome this difficulty, we could try to apply a trick of Schoen-Yau \cite{sy1} 
to compactify these ends, however it does not work in our situation. Indeed, the first step in Schoen-Yau's method is to deform the metric $\what{g}$
on the parts which are asymptotic to the cylinders so that these components {\it coincide} with the cylinders. However this transformation destroys Equation (\ref{schoenyau}) 
and this adds a real difficulty for our purpose. Nevertheless, it seems reasonable to think that our result is true in this more general context.


\subsection{Proof of Theorem \ref{spacetime}}


{}From the work \cite{yau1} and since we assumed that $\Omega$ has no apparent horizon in its interior, there exists a solution $u$ to the Jang Equation 
(\ref{jang}) defined and smooth on $\Omega$ with Dirichlet boundary condition 
\begin{eqnarray*}
 u_{|\Sigma}\equiv 0.
\end{eqnarray*}
This boundary condition ensures that the metric $\what{g}$ coincides with the metric $g$ on the boundary $\Sigma$ so that the Dirac operators 
$D^\SSi$ acting on $\SSi$ and $D^{\what{\SSi}}$ on $\what{\SSi}$ also coincide. In particular, we have $\lambda_1(D^{\what{\SSi}})=\lambda_1(D^\SSi)=\lambda_1(D^\Sigma)$ 
where the last equality comes from the fact that $\Sigma$ is a $2$-dimensional manifold and then $D^\SSi=D^\Sigma$. Moreover, from a calculation in \cite{yau1}
we have:
\begin{eqnarray*}
\what{H}_{\Sigma}-\what{g}(X,\what{N})=f^{-1}H_{\Sigma}-\sigma|\nabla u|{\rm Tr}_\Sigma(K) 
\end{eqnarray*}
where $\sigma\in\{\pm 1\}$ and thus 
\begin{eqnarray}\label{boundblack}
\what{H}_{\Sigma}-\what{g}(X,\what{N})\geq \sqrt{H_\Sigma^2-{\rm Tr}_\Sigma(K)^2}
\end{eqnarray}
holds since $f=-\<\partial_t,\what{\nu}\>=1/\sqrt{1+|\nabla u|^2}$. Here $\what{N}$ denotes the unit outward normal vector field of 
$\Sigma$ in $\what{\Omega}$. On the other hand, the resulting Riemannian manifold $\what{\Omega}$ satisfies  condition 
(\ref{scalarcond}) because of (\ref{schoenyau}) and where the vector field $X$ is defined by (\ref{defX}). The assumptions of Theorem \ref{newextrin} are 
fulfilled and then
\begin{eqnarray}\label{new1}
\lambda_1(D^\Sigma)\geq\frac{1}{2}\underset{\Sigma}{\inf}\big(\what{H}_{\Sigma}-\what{g}(X,\what{N})\big).
\end{eqnarray} 
Combining (\ref{boundblack}) and (\ref{new1}) give the estimate.

\medskip

Now assume that equality is achieved. Once again we apply Theorem \ref{newextrin} and then $\what{\Omega}$ has a parallel spinor field $\Phi$. In particular, 
$\what{\Omega}$ is Ricci flat and since it is a $3$-dimensional domain, it is flat. Moreover, if we have equality in (\ref{schoenyau}), then the second 
fundamental form of $\what{\Omega}$ in $M\times\mathbb{R}$ is $K_{ij}$. So we can choose a  coordinates system  $\what{x}=(\what{x}_1,\what{x}_2,\what{x}_3)$ 
in a neighborhood $\mathcal{U}$ of a point $p\in\Omega$ such that $\what{g}_{ij}=\delta_{ij}$. In this chart, we have:
\begin{eqnarray*}
g_{ij}=\delta_{ij}-\frac{\partial u}{\partial\what{x}_i}\frac{\partial u}{\partial\what{x}_j}
\end{eqnarray*}

and this shows that if $(\what{x}_1,\what{x_2},\what{x_3},t)$ are coordinates in the Minkowski spacetime
\begin{eqnarray*}
\mathbb{R}^{3,1}=\big(\mathbb{R}^4,\sum_{i=1}^3d\what{x_i}^2-dt^2\big)
\end{eqnarray*}

the graph of $u$ over $\mathcal{U}$ isometrically embeds in $\mathbb{R}^{3,1}$ with second fundamental form given by $K_{ij}$. 
\hfill$\square$

\begin{remark}\label{modlocmass}
This result leads to the definition of a new quasi-local mass for $2$-spheres with positive Gauss curvature similar to Liu-Yau \cite{liuyau}. The idea
here, following \cite{zhang} and \cite{rollmops}, consists in localizing Witten's proof of the positive mass theorem \cite{witten}. Namely let $\Sigma$ denotes a $2$-sphere with positive Gauss curvature which bounds a compact domain $\Omega$ in an initial data set $(M^3,g,K)$. Since
$\Sigma$ has positive Gauss curvature, the Weyl's embedding theorem (\cite{nirenberg} or \cite{pogorelov}) ensures the existence of an isometric embedding 
(unique up to the isometries of $\mathbb{R}^3$) of $(\Sigma^2,g)$ in the flat Euclidean space with positive mean curvature $H^0_\Sigma$. Then, using Formula 
(\ref{dira-extr2}), the restriction of a parallel spinor field $\Phi\in\Gamma\big(\mathbb{S}(\mathbb{R}^3)\big)$ to $\Sigma$ gives a solution to
the boundary Dirac equation 
\begin{eqnarray*}
D^\Sigma\Phi=\frac{1}{2}H^0_\Sigma\Phi.
\end{eqnarray*}

We then define the quantity 
\begin{eqnarray*}
\mathcal{Q}(\Sigma,\Psi):=\int_{\Sigma}H^0_\Sigma|\Psi|^2ds-\int_{\Sigma}\sqrt{H^2_\Sigma-\mathrm{Tr}_{\Sigma}(h)^2}|\Psi|^2ds
\end{eqnarray*}

where $\Psi\in\Gamma\big({\mathbb S}(\what{\Omega})\big)$ is a solution of the boundary problem
$$\left\lbrace
\begin{array}{ll}
\what{D}\Psi = 0 & \qquad\text{on}\,\,\what\Omega\\
\what{B}^{\pm}\Psi_{|\Sigma}= \what{B}^{\pm}\Phi& \qquad\text{along}\,\,\what{\Sigma}
\end{array}
\right.$$

on the graph $\what{\Omega}$ whose existence is given by Lemma \ref{existenceMIT}. Then, following the proof of Theorem \ref{spacetime} and using 
Theorem \ref{FunctEst}, we obtain:
\begin{eqnarray}\label{posquas}
\mathcal{Q}(\Sigma,\Psi)\geq 0.
\end{eqnarray}

The new quasi-local mass is then defined by
\begin{eqnarray*}
 m(\Sigma):=\min_{\mathcal{P}}\mathcal{Q}(\Sigma,\Psi)
\end{eqnarray*}

where $\mathcal{P}$ denotes the space of constant spinors with unit norm. Form (\ref{posquas}), it is clear that $m(\Sigma)\geq 0$ and if $m(\Sigma)=0$, 
it is straightforward to see, from the proof of Theorem \ref{spacetime}, that $\Sigma$ embeds in the Minkowski spacetime. 
\end{remark}


\subsection{Rigidity for the constraint equations}


We provide here a direct application of our main estimate to get a rigidity result for the constraint equation with spherical boundary. More precisely, 
we will say that a compact domain $\Omega$ of an initial data set $(M^3,g,K)$ has a spherical boundary if the boundary of $\Omega$ endowed with the 
induced metric is isometric to a round sphere. In the following, we can assume without loss of generality, that the boundary sphere has radius one. 
The two following domains:
\begin{enumerate}
\item the unit Euclidean ball $(\mathcal{B}^3,eucl)$ and 
\item the spherical cap defined by 
\begin{eqnarray*}
\mathcal{C}_r^3:=\{x=(x_1,x_2,x_3,t)\in\mathbb{R}^{4}\,/\,x^2_1+x^2_2+x^2_3-t^2=-r^2,\,0<t\leq\sqrt{1+r^2}\} 
\end{eqnarray*}
\end{enumerate}

give examples of compact domains with spherical boundary. It is interesting to point out that these two domains can be seen as domains in spacelike hypersurfaces of 
the Minkowski spacetime: a totally geodesic hyperplane $\mathbb{R}^3$ for the unit Euclidean ball and a totally umbilical hyperbolic space 
$\mathbb{H}^3(-\frac{1}{r^2})$ with constant curvature equals to $-(1/r^2)$ for the spherical caps. Moreover, it is clear that, in both cases, 
the mean curvature vector has constant length equal to 2. In fact, this comes from a more general fact:
\begin{theorem}\label{uniqueness}
Let $\Omega$ be a compact domain with spherical boundary in an initial data $(M^3,g,K)$ and assume that $\Omega$ has no apparent horizon in
its interior. If the mean curvature vector of $\Sigma$ is such that $|\mathcal{H}|\geq 2$, then $\Omega$ embeds as a spacelike hypersurface in the 
Minkowski spacetime. 
\end{theorem}

{\it Proof:}
First, note that since we assume that the boundary is spherical we have $\lambda_1(D^\Sigma)=1$ and the associated eigenspinor is a real Killing spinor. On 
the other hand, we can apply Theorem \ref{spacetime} and the assumption on the norm of the mean curvature vector to conclude that 
$\lambda_1(D^\Sigma)\geq 1$ and so equality occurs in (\ref{eigenspacetime}). From the equality case, we have that $\Omega$ can be locally embedded in
the Minkowski spacetime. As seen in the proof of Theorem \ref{spacetime}, this embedding is constructed from a solution of the Jang equation 
$u\in C^\infty(\Omega)$. Here we easily see that the manifold $\what{\Omega}=(\Omega,\what{g})$, with $\what{g}_{ij}=g_{ij}+u_iu_j$, is a compact flat Riemannian manifold with spherical boundary and
so it is isometric to the unit Euclidean ball. In particular, $\Omega$ is simply connected and thus $(\Omega,g)$ can be globally embedded in the Minkowski 
spacetime. 
\hfill$\square$\\

If we assume that the extrinsic mean curvature (that is ${\rm Tr}(K)$) is constant, we prove that the only domains with spherical boundary are the Euclidean 
balls and the spherical caps. In particular, this statement provides an unified approach to previous results of Miao \cite{miao}, Shi-Tam \cite{shitam}, 
\cite{shitam1} and the author \cite{rollmops} for the rigidity of Euclidean balls and spherical caps. More precisely, we prove:
\begin{corollary}\label{uniqueness1}
Under the assumptions of Theorem \ref{uniqueness} and if the extrinsic mean curvature is 
constant, then $\Omega$ is isometric to an Euclidean ball (if ${\rm Tr}(K)=0$) or to a spherical cap (if ${\rm Tr}(K)\neq 0$).
\end{corollary}

{\it Proof:}
{}From Theorem \ref{uniqueness} and since the extrinsic mean curvature is constant, the domain $\Omega$ is a constant mean curvature
spacelike hypersurface with spherical boundary in the Minkowski spacetime. By a result of Alias and Pastor \cite{alias1}, $\Omega$ is 
an Euclidean ball or a spherical cap. 
\hfill$\square$\\

Using a result of Alias and Malacarne \cite{alias2} similar to \cite{alias1}, it is straightforward to see that the statement of Corollary 
\ref{uniqueness1} is also true for initial data sets with constant higher order mean curvature. In particular, it holds if the scalar curvature of 
$\Omega$ is constant.


\subsection{Surfaces in Minkowski spacetime}


In this section, we prove that Inequality (\ref{eigenspacetime}) holds in the case of surfaces embedded in the Minkowski spacetime with spacelike mean 
curvature vector without any assumption on the existence of apparent horizon. Moreover, since this estimate only involves  geometrical data of the surface $\Sigma$ and does not depend on the 
bounding domain $\Omega$, we focus especially on the geometrical properties of the surfaces for which (\ref{eigenspacetime}) is an equality. More 
precisely,  following \cite{miaoshitam}, we prove
\begin{theorem}\label{Mink}
Let $\Sigma$ be a closed, connected surface with spacelike mean curvature vector field which bounds a compact domain of a spacelike hypersurface in the 
Minkowski spacetime $\mathbb{R}^{3,1}$, then Inequality (\ref{eigenspacetime}) holds. Moreover, equality occurs if and only if $\Sigma$ is a round sphere 
in $\mathbb{R}^{3,1}$.
\end{theorem}

Before giving the proof of this result, we fix some notations and recall some basic facts on spinors in the Minkowski spacetime. In the following, 
the triplet $(\SM,\nabla^{\RM},\widetilde{\gamma})$ represents the Dirac bundle over $\RM$ made of the complex spinor bundle $\SM$, the spin Levi-Civita connection $\nabla^{\RM}$ 
acting on $T\RM\times\SM$ and the Clifford multiplication $\widetilde{\gamma}$ acting on $T\RM\otimes\SM$. If $\Omega$ denotes a $3$-dimensional spacelike 
hypersurface of $\RM$ with second fundamental form $K$, we define the {\it hypersurface} Dirac bundle $(\SM_{|\Omega},\nabla,\gamma)$ 
on which the following spinorial Gauss formula holds:
\begin{eqnarray}\label{spingaussformule}
\nabla^{\RM}_X\varphi=\nabla_X\varphi-\frac{1}{2}\widetilde{\gamma}(KX)\widetilde{\gamma}(T)\varphi.
\end{eqnarray}

Here $T$ is the unit timelike vector field normal to $\Omega$. It is not difficult to see that the Clifford multiplications $\widetilde{\gamma}$ and $\gamma$
are related by
\begin{eqnarray}\label{idclif}
\gamma(X)\varphi=i\widetilde{\gamma}(X)\widetilde{\gamma}(T)\varphi
\end{eqnarray}

for all $X\in\Gamma(T\Omega)$ and $\varphi\in\Gamma(\SM_{|\Omega})$. We should note that $\Omega$ is also endowed with an {\it intrinsic} Dirac bundle
$(\mathbb{S}(\Omega),\nabla^\Omega,\gamma^\Omega)$ satisfying:
\begin{eqnarray*}
\big(\SM_{|\Omega},\nabla,\gamma\big)\simeq\big(\mathbb{S}(\Omega)\oplus\mathbb{S}(\Omega),\nabla^\Omega\oplus\nabla^\Omega,
\gamma^\Omega\oplus-\gamma^\Omega\big)
\end{eqnarray*}

since $\RM$ is four-dimensional. To simplify our notations, we identify these two bundles and only work with the extrinsic one. Similarly, the surface $\Sigma$ is endowed with
two Dirac bundles, an intrinsic one $({\bf S}:=\mathbb{S}\Omega_{|\Sigma},\nabla^{{\bf S}},\gamma^{\bf S})$ and an extrinsic one 
$(\mathbb{S}\Sigma,\nabla^\Sigma,\gamma^\Sigma)$. In this situation, since $\Sigma$ is $2$-dimensional, we have the following identification
\begin{eqnarray*}
({\bf S}:=\mathbb{S}\Omega_{|\Sigma},\nabla^{{\bf S}},\gamma^{\bf S})\simeq(\mathbb{S}\Sigma,\nabla^\Sigma,\gamma^\Sigma).
\end{eqnarray*}
Moreover, we also have the following identification for the Clifford multiplications:
\begin{eqnarray}
 \gamma^{\SSi}(X)\varphi=\gamma(X)\gamma(N)\varphi=\widetilde{\gamma}(X)\widetilde{\gamma}(N)\varphi
\end{eqnarray}
for all $X\in\Gamma(T\Sigma)$ and $\varphi\in\Gamma(\mathbf{S})$. Combining the spin Gauss formula of the immersions $\Sigma\hookrightarrow\Omega$ and
$\Omega\hookrightarrow\RM$ give that, for all $X\in\Gamma(T\Sigma)$ and $\varphi\in\Gamma(\SSi)$, we have
\begin{eqnarray*}
 \nabla^{\RM}_X\varphi=\nabla^\SSi_X\varphi+\frac{1}{2}\widetilde{\gamma}(AX)\widetilde{\gamma}(N)\varphi-\frac{1}{2}\widetilde{\gamma}(KX)\widetilde{\gamma}
(T)\varphi.
\end{eqnarray*}
Using this identity, it is straightforward to check that the Dirac operator $D^\SSi$ satisfies
\begin{eqnarray}\label{DirMin}
 D^\SSi\varphi=\frac{1}{2}H_\Sigma\varphi+\sum_{j=1}^2\widetilde{\gamma}(e_j)\widetilde{\gamma}(N)\nabla^{\RM}_{e_j}\varphi+
\frac{1}{2}\sum_{j=1}^2\widetilde{\gamma}(e_j)\widetilde{\gamma}(N)\widetilde{\gamma}(Ke_j)\widetilde{\gamma}(T)\varphi
\end{eqnarray}
where $\{e_1,e_2\}$ is a local $g$-orthonormal frame of $T\Sigma$. 

\vspace{0.5cm}

{\it Proof of Theorem \ref{Mink}:} 
{}From Lemma $4.1$ in \cite{miaoshitam}, $\Sigma$ spans a compact, smoothly immersed, maximal hypersurface $\Omega$ in $\mathbb{R}^{3,1}$. Denote 
by $K$ the second fundamental form of $\Omega$ in $\mathbb{R}^{3,1}$ (which satisfies
${\rm Tr}(K)=0$ since $\Omega$ is maximal), the Gauss formula gives $R=|K|^2\geq 0$. Here
$R$ is the scalar curvature of $\Omega$ equipped with the metric induced by the Minkowski spacetime. Moreover, the 
mean curvature vector of $\Sigma$ satisfies
\begin{eqnarray}\label{mcv}
 |\mathcal{H}|^2=H_\Sigma^2-{\rm Tr}_{\Sigma}(K)^2\leq H_\Sigma^2.
\end{eqnarray}

{}From Lemma $4.2$ in \cite{miaoshitam} and since $\Sigma$ has a spacelike mean curvature vector, we have $H_\Sigma>0$. Now since $R\geq 0$ and $H_\Sigma>0$,
we can apply the Hijazi-Montiel-Zhang estimate \cite{hijazi.montiel.zhang:1} (which is inequality (\ref{ne1}) with $X\equiv 0$) to get
\begin{eqnarray*}
 \lambda_1(D^\Sigma)\geq\frac{1}{2}\inf_\Sigma(H_\Sigma)\geq\frac{1}{2}\inf_{\Sigma}|\mathcal{H}|
\end{eqnarray*}

because of (\ref{mcv}). Assume now that equality is achieved. From the equality case of (\ref{ne1}), we deduce that $\Omega$ has a parallel spinor
(in particular it is flat since Ricci flat and $3$-dimensional) and that the mean curvature $H_\Sigma$ is constant. Moreover since 
$R=|K|^2=0$, the domain $\Omega$ is totally geodesic in $\mathbb{R}^{3,1}$. In particular, it is straightforward to see that
the value of the mean curvature $H_\Sigma$ is $2\lambda_1(D^{\Sigma})$. We also know that the eigenspace associated with the eigenvalue $\lambda_1(D^\Sigma)$ is 
given by the restriction to $\Sigma$ of the space of parallel spinor fields on $\Omega$. Let $\Phi_0\in\Gamma(\SM_{|\Omega})$ such a 
parallel spinor and define for $p\in\Sigma$:
\begin{eqnarray*}
 \Psi_0(p):=\widetilde{\gamma}\big(\frac{1}{2}H_{\Sigma}\xi(p)+N(p)\big)\Phi_0.
\end{eqnarray*}

Here $\xi$ denotes the position vector field in $\RM$. Using Identity (\ref{DirMin}) and the fact that $\Omega$ is totally geodesic in $\RM$, 
we compute that:
\begin{eqnarray*}
 D^\SSi\Psi_0= \frac{1}{2}H_\Sigma\Psi_0+{\rm Tr}(\frac{H_\Sigma}{2}\Id-A)\Phi_0=\frac{1}{2}H_\Sigma\Psi_0
\end{eqnarray*}

and thus $\Psi_0$ is an eigenspinor for the Dirac operator $D^\SSi$ associated with $\lambda_1(D^\Sigma)$. Then $\Psi_0$ extends to a parallel spinor on $\Omega$ and we 
compute for $X\in\Gamma(T\Sigma)$:
 \begin{eqnarray*}
  0=\nabla_X\Psi_0 & = & \nabla^{\RM}_X\big(\widetilde{\gamma}\big(\frac{1}{2}H_{\Sigma}\xi(p)+N(p)\big)\Phi_0\big)\\
& = & \widetilde{\gamma}\big(\frac{H_\Sigma}{2}X-AX\big)\Phi_0.
 \end{eqnarray*}

Since $\Phi_0$ has no zeros, we immediately get $AX=(H_\Sigma/2)X$ for all $X\in\Gamma(T\Sigma)$. From the Gauss formula of the immersion $\Sigma\hookrightarrow\RM$,
we observe that the Gauss curvature of $\Sigma$ is constant equal to
\begin{eqnarray*}
 2K_\Sigma=|\mathcal{H}|^2-|II|^2=\frac{H_\Sigma^2}{2}>0
\end{eqnarray*}
and so $\Sigma$ is isometric to a round sphere.

\medskip

If now we assume that $\Sigma$ is a round sphere with radius $r>0$, it embeds in a totally geodesic hyperplane $\mathbb{R}^3$ such that $|\mathcal{H}|^2=\frac{4}{r^2}$.
On the other hand, the first eigenvalue of the Dirac operator is $\frac{1}{r}$ and so equality is achieved in (\ref{eigenspacetime}).
\hfill$\square$\\

For our geometric applications, we prove a new upper bound for this first eigenvalue in terms of geometrical data of $\Sigma$ stated as follow:
\begin{proposition}\label{upper}
Let $\Sigma$ be a compact and connected spacelike surface isometrically immersed in the Minkowski spacetime with spacelike mean curvature vector. Assume that 
$\Sigma$ is contained in a spacelike hypersurface $M$ of $\mathbb{R}^{3,1}$ such that the trace of its second fundamental form $K$ along $\Sigma$ is constant.
Then the first eigenvalue of the Dirac operator $D^\Sigma$ satisfies
\begin{eqnarray}
\lambda_1(D^\Sigma)^2\leq\frac{1}{4}\sup_{\Sigma}\big(H^2_\Sigma-{\rm Tr}_{\Sigma}(K)^2\big). 
\end{eqnarray}

Moreover, if equality occurs, $(1/2)|\mathcal{H}|=(1/2)\sqrt{H^2_\Sigma-{\rm Tr}_{\Sigma}(K)^2}$ is the first eigenvalue of the Dirac operator of $\Sigma$.
\end{proposition}

\begin{remark}
It is straightforward to check that this result holds for codimension two compact submanifolds isometrically immersed in a spacelike hypersurface of the 
Minkowski space $\mathbb{R}^{n,1}$. In this situation, our estimate generalizes a previous result of Ginoux \cite{ginoux} for hypersurfaces in the 
Hyperbolic space.
\end{remark}

{\it Proof of Proposition \ref{upper}:}
We consider the Dirac-Witten type operator defined by
\begin{eqnarray*}
D^W\varphi:=\sum_{j=1}^2\gamma(e_j)\nabla_{e_j}\varphi.
\end{eqnarray*}

for $\varphi\in\Gamma(\mathbf{S})$. By a direct calculation, we have:
\begin{eqnarray*}
 D^W\varphi=\gamma(N)\Big(D^{\mathbf{S}}\varphi-\frac{H_\Sigma}{2}\varphi\Big)
\end{eqnarray*}
where $N$ is the unit inner vector field normal to $\Sigma$ in $M$. Moreover, since
\begin{eqnarray*}
 D^W\big(\gamma(N)\varphi\big)=H_\Sigma\varphi-\gamma(N)D^W\varphi
\end{eqnarray*}
we compute
\begin{eqnarray*}
 (D^{\mathbf{S}})^2\varphi=(D^W)^2\varphi+\frac{1}{2}\gamma(d^\Sigma H_\Sigma)\gamma(N)\varphi+\frac{H^2_\Sigma}{4}\varphi
\end{eqnarray*}
and the Rayleigh quotient for the first eigenvalue $\lambda_1(D^{\mathbf{S}})^2$ writes
\begin{equation}\label{rayl}
 \lambda_1(D^{\bf{S}})^2\leq\frac{\int_\Sigma\<(D^{\mathbf{S}})^2\varphi,\varphi\>ds}{\int_\Sigma|\varphi|^2ds}
=\frac{\int_\Sigma\<(D^W)^2\varphi+\frac{H^2_\Sigma}{4}\varphi,\varphi\>ds}{\int_\Sigma|\varphi|^2ds}
\end{equation}
for all $\varphi\in\Gamma(\mathbf{S})$ since $\<\gamma(d^\Sigma H_\Sigma)\gamma(N)\varphi,\varphi\>$ is purely imaginary. Now, since $M$ is an immersed
spacelike hypersurface of $\mathbb{R}^{3,1}$, we consider the restriction to $M$ of a constant spinor field $\Phi_0\in\Gamma(\SM)$ 
which satisfies by the spinorial Gauss formula (\ref{spingaussformule})
\begin{eqnarray*}
 \nabla_X\Phi_0=\frac{1}{2}\widetilde{\gamma}(KX)\widetilde{\gamma}(T)\Phi_0
\end{eqnarray*}
for all $X\in\Gamma(TM)$ and where $T$ denotes the unit timelike future oriented vector field normal to $M$. From the identification (\ref{idclif}) of 
the Clifford multiplications $\gamma$ and $\widetilde{\gamma}$ we have:
\begin{eqnarray*}
 \nabla_X\Phi_0=-\frac{i}{2}\gamma(KX)\Phi_0
\end{eqnarray*}
for all $X\in\Gamma(TM)$. Now we compute
\begin{eqnarray*}
 D^W\Phi_0=\frac{i}{2}{\rm Tr}_\Sigma(K)\Phi_0
\end{eqnarray*}

and since ${\rm Tr}_{\Sigma}(K)$ is constant, we conclude that 
\begin{eqnarray*}
 (D^W)^2\Phi_0 = -\frac{{\rm Tr}_\Sigma(K)^2}{4}\Phi_0.
\end{eqnarray*}

Using $\Phi_0$ as a test-spinor in (\ref{rayl}) gives immediately the estimate. If equality is achieved, it is straightforward to see that the norm of the 
mean curvature vector is the first eigenvalue of $D^{\mathbf{S}}$.
\hfill$\square$\\

As an application of our previous estimates, we obtain the following uniqueness result:
\begin{corollary}\label{genalex}
Let $\Sigma$ be a closed, connected, spacelike surface in the Minkowski spacetime. If $\Sigma$ bounds a compact spacelike domain $\Omega$ 
such that $H_\Sigma$ and ${\rm Tr}_\Sigma(K)$ are constant, then $\Sigma$ is a round sphere.
\end{corollary}

{\it Proof:}
{}From our assumptions, we can apply Theorem \ref{Mink} and  Proposition \ref{upper} to get:
\begin{eqnarray*}
 \frac{1}{4}\inf_{\Sigma}\big(H^2_\Sigma-{\rm Tr}_{\Sigma}(K)^2\big)\leq\lambda_1(D^\Sigma)^2\leq\frac{1}{4}\sup_{\Sigma}\big(H^2_\Sigma-{\rm Tr}_{\Sigma}(K)^2\big)
\end{eqnarray*}

On the other hand, since $H_\Sigma$ and ${\rm Tr}_\Sigma(K)$ are constant, the norm of the mean curvature vector is also constant so that we have 
equality in the previous inequalities. The rigidity statement of Theorem \ref{Mink} allows to conclude.
\hfill$\square$\\

This result gives a unified proof of the Euclidean and Hyperbolic version of the Alexandrov Theorem for surfaces (see 
\cite{hijazi.montiel.zhang:1} and \cite{hijazi.montiel.roldan:02} for example). Namely, we have:
\begin{corollary}\label{alex}
The only compact and connected surfaces with constant mean curvature in the Euclidean and the Hyperbolic space are the round spheres. 
\end{corollary}

{\it Proof:}
It is enough to note that the Euclidean space embeds as a totally geodesic spacelike hyperplane in the Minkowski spacetime and in the Hyperbolic space as a totally umbilical with 
constant mean curvature spacelike hypersurface. It is then clear that in both cases, the surface $\Sigma$ satisfies the assumptions of  Corollary \ref{genalex} so that 
we get immediately the result.
\hfill$\square$\\


\bibliographystyle{amsalpha}


\end{document}